\documentclass[12pt]{article}
\usepackage{latexsym}

\usepackage[tbtags]{amsmath}
\usepackage{epsfig,amstext,amssymb,amsthm,latexsym}
\usepackage{amssymb,color}
\usepackage{tikz}
\usetikzlibrary{calc,spy}
\usepackage{pgfplots}
\UseRawInputEncoding

\definecolor{c20}{rgb}{0.,0.7,0.}
\definecolor{c30}{rgb}{0.,0.,1.}
\definecolor{c40}{rgb}{1,0.1,0.7}
\definecolor{c50}{rgb}{1,0,0}

\date{}
\newtheorem{lemma}{Lemma}[section]

\makeatletter 
\@addtoreset{equation}{section}
\makeatother 

\begin{document}
\title{ Note on the First Part of Borel-Cantelli Lemma}

\author{Alexei Stepanov  \thanks{\noindent Higher School of Computer Science and Applied Mathematics,  Education and Research Cluster "Institute of High Technology``,\  Immanuel Kant Baltic Federal University,  A. Nevskogo  14, Kaliningrad, 236041 Russia; email: alexeistep45@mail.ru}}

\maketitle
\begin{abstract} 
In this short note, we briefly discuss the Borel-Cantelli lemma and propose a new generalization of the first part of  it. 
\end{abstract}

\noindent {\it MSC2020-Mathematical Sciences Classification System} 60, 62.

\section{ Introduction}
Suppose $A_1,A_2,\ldots$ is a sequence of events on a common probability space and that $A^c_i$ denotes the complement of event $A_i$. The Borel-Cantelli lemma, presented below as Lemma~\ref{lemma1.1}, is used  for producing strong limit results.
\begin{lemma}\label{lemma1.1}
\begin{enumerate}
\item If, for any sequence  $A_1,A_2,\ldots$ of events,
\begin{equation}\label{1.1}
\sum_{n=1}^\infty P(A_n)<\infty,
\end{equation}
then $P(A_n\ i.o.)=0$;
 \item If $A_1,A_2,\ldots$ is a sequence
of independent events and if $\sum_{n=1}^\infty P(A_n)=\infty$,
then $P(A_n\ i.o.)=1$.
\end{enumerate}
\end{lemma}
The independence condition in the second part of Lemma~\ref{lemma1.1} has been weakened by a number of authors, including
Chung and Erdos (1952), Erdos and Renyi (1959), Lamperti (1963), Kochen and Stone (1964), Spitzer (1964), Chandra (1999, 2008),  Petrov (2002, 2004) and others.

The first part of  Borel-Cantelli lemma has been generalized in Barndorff-Nielsen (1961) and  Balakrishnan and Stepanov (2010). For
a review on the Borel–-Cantelli lemma, one may refer to the book of Chandra (2012). The result of Barndorff-Nielsen is presented below as Lemma~\ref{lemma1.2}.
\begin{lemma}\label{lemma1.2}
Let $A_n\ (n\geq 1)$ be a sequence of events such that $P(A_n)\rightarrow 0$. If
\begin{equation}\label{1.2}
\sum_{n=1}^\infty P(A_n A^c_{n+1})<\infty,
\end{equation}
then $P(A_n\ i.o.)=0$.
\end{lemma}
Observe  that condition (\ref{1.2})  in Lemma~\ref{lemma1.2} is weaker than  condition (\ref{1.1}) in the Borel-Cantelli lemma. It should be noted that  Lemma~\ref{lemma1.2} was proposed in    Barndorff-Nielsen (1961) with the wrong proof. For explanations and the correct proof of Lemma~\ref{lemma1.2}, see Balakrishnan and Stepanov (2021) and Stepanov and Dembinska (2022). We also present here the result of Balakrishnan and Stepanov (2010). 
\begin{lemma}\label{lemma1.3}
Let  $A_1,A_2,\ldots$ be a sequence of events such that $P(A_n)\rightarrow 0$.  If, for some $m\geq0$,
$$
\sum_{n=1}^\infty P(A^c_nA_{n+1}^c\ldots A_{n+m-1}^c A_{n+m})<\infty,
$$
then $P(A_n\ i.o.)=0$.
\end{lemma}

\section{New Result}
The following result holds true.
\begin{lemma}\label{lemma2.1}
Let  $A_1,A_2,\ldots$ be a sequence of events such that $P(A_n)\rightarrow 0$.  If, for some $m\geq0$,
$$
\sum_{n=1}^\infty P(A_nA_{n+1}^c\ldots A_{n+m-1}^c A_{n+m}^c)<\infty,
$$
then $P(A_n\ i.o.)=0$.
\end{lemma}

\begin{gproof}{} For  a large fixed $k$ and some $n\geq 1$, we have
\begin{eqnarray*}
&&\hspace{-5ex}P(A_n+A_{n+1}+\ldots+A_{n+k})=P(A_n A^c_{n+1}\ldots A_{n+k}^c)+P(A_{n+1}+\ldots+A_{n+k})\\
&=& P(A_nA_{n+1}^c\ldots A^c_{n+k})+P(A_{n+1}A_{n+2}^c\ldots A_{n+k}^c)+P(A_{n+2}+\ldots+A_{n+k})\\
&=&\ldots\\
&=&P(A_n A^c_{n+1}\ldots A_{n+k}^c)+P(A_{n+1} A^c_{n+2}\ldots A_{n+k}^c)+\ldots+P(A_{n+k-1}A_{n+k}^c)+ P(A_{n+k}).
\end{eqnarray*}
Then, for a small $m\geq 0$, we get that
\begin{eqnarray*}
&&\hspace{-5ex}P(A_n+A_{n+1}+\ldots+A_{n+k})\\
&\leq& P(A_n A^c_{n+1}\ldots A_{n+m}^c)+P(A_{n+1} A^c_{n+2}\ldots A_{n+m+1}^c)+\ldots+ P(A_{n+k-m}A_{n+k-m+1}^c\ldots A_{n+k}^c)+o_n(1).
\end{eqnarray*}
Letting $k\rightarrow \infty$, we obtain
$$
P\left(\sum_{i=n}^\infty A_i\right)\leq o_n(1)+\sum_{i=n}^\infty P(A_i A^c_{i+1}\ldots A_{i+m}^c).
$$ 
The result readily follows.
\end{gproof}
Observe that Lemma~\ref{lemma2.1} generalizes Lemma~\ref{lemma1.2}. 

\section*{References}

\begin{description} 
\item Balakrishnan, N., Stepanov, A. (2010).\ Generalization of Borel-Cantelli lemma. {\it The Mathematical Scientist}, {\bf 35}, 61--62.

\item Balakrishnan, N., Stepanov, A. (2021).\ A note on the Borel–-Cantelli lemma. arXiv:2112.00741v2, [math.PR].

\item Barndorff-Nielsen, O. (1961).\ On the rate of growth of the partial maxima of a sequence of independent identically distributed random variables. {\it Math. Scand.},  9, 383--394.

\item Chandra, T.K. (1999).\ {\it A First Course in Asymptotic Theory of Statistics}. Narosa Publishing House Pvt. Ltd., New Delhi.

\item Chandra, T.K.  (2008).\ Borel-Cantelli lemma under dependence conditions. {\it Statist. Probab. Lett.}, {\bf 78}, 390–-395.

\item Chandra, T.K., (2012).\ {\it The Borel–-Cantelli Lemma}. Springer Briefs in Statistics.

\item Chung, K.L. and Erdos, P. (1952).\ On the application of the Borel-Cantelli lemma. {\it Trans.  Amer. Math. Soc.},  72,
179--186.

\item Erdos, P. and Renyi, A. (1959).\ On Cantor's series with convergent $\sum 1/q_n$. {\it Ann. Univ. Sci. Budapest. Sect.
Math.},  2, 93--109.

\item Kochen, S.B. and Stone, C.J. (1964).\ A note on the Borel-Cantelli lemma. {\it Illinois J. Math.},  8, 248--251.

\item Lamperti, J. (1963).\ Wiener's test and Markov chains. {\it J. Math. Anal. Appl.},  6, 58--66.

\item Petrov, V.V. (2002).\ A note on the Borel-Cantelli lemma. {\it Statist. Probab. Lett.}, {\bf 58}, 283--286.

\item  Petrov, V.V. (2004).\ A generalization of the Borel--Cantelli lemma. {\it Statist. Probab. Lett.}, {\bf 67}, 233--239.

\item Spitzer, F. (1964).\ {\it Principles of Random Walk}.  Van Nostrand, Princeton, New Jersey.

\item Stepanov A. (2014).\ On the Use of the Borel-Cantelli Lemma   in Markov Chains, {\it Statist. Probab. Lett.}, {\bf 90}, 149--154.

\item Stepanov, A.,  Dembinska, A. (2023).\ Limit Theorems for  the  Uppermost $m$-th Spacing Based on Weak Geometric   Records,
 {\it Statist. Probab. Lett.}, published online 2022, doi.org/10.1016/j.spl.2021.109351.

\end{description}

\end{document}